\documentclass[12pt]{article}
\usepackage{amsmath}
\usepackage{amssymb}
\numberwithin{equation}{section}
\usepackage{epic}

\title{A PDE for the joint distributions of \\ the Airy Process}

\author{
M. Adler\thanks{ Department of Mathematics, Brandeis
University, Waltham, Mass 02454, USA. E-mail:
adler@brandeis.edu.  The support of a National Science
Foundation grant \# DMS-01-00782 is gratefully
acknowledged.}~~~~~~ P. van Moerbeke\thanks{ Department
of Mathematics, Universit\'e de Louvain, 1348
Louvain-la-Neuve, Belgium and Clay Mathematics
Institute, One Bow Street, Cambridge, MA 02138, USA.
E-mail: vanmoerbeke@math.ucl.ac.be . The support of a
National Science Foundation grant \# DMS-01-00782,
European Science Foundation, Nato, FNRS and Francqui
Foundation grants is gratefully acknowledged.}}

\date{}

\newcommand{\MAT}[1]{\left(\begin{array}{*#1c}}
\newcommand{\mat}{\end{array}\right)}
\newcommand{\qed}{\leavevmode\unskip\nobreak\penalty200\hskip2pt\null
\nobreak\hfill\rule{1.1ex}{1.1ex}
\medbreak }

\newcommand{\rg}{\rightarrow}

\newcommand{\LR}{{\cal L}}

\newcommand{\iy}{\infty}
\newcommand{\pl}{\partial}
\newcommand{\al}{\alpha}

\newcommand{\no}{\nonumber}

\newenvironment{proof}{\medskip\noindent{\it Proof:\/} }{\qed}
\newenvironment
        {remark}{\medskip\noindent\underline{\it Remark:\/} }{\medbreak}

\newcommand{\Dt}{\Delta}

\newcommand{\BR}{{\mathbb R}}
\newcommand{\lb}{\lambda}

\def\be#1\ee{\begin{equation}#1\end{equation}}
\def\bea#1\eea{\begin{eqnarray}#1\end{eqnarray}}
\def\bean#1\eean{\begin{eqnarray*}#1\end{eqnarray*}}

\newcommand{\Tr}{\operatorname{\rm Tr}}

\newtheorem{definition}{Definition}[section]
\newtheorem{theorem}[definition]{Theorem}
\newtheorem{lemma}[definition]{Lemma}
\newtheorem{corollary}[definition]{Corollary}

\catcode `!=11

\newdimen\squaresize
\newdimen\thickness
\newdimen\Thickness
\newdimen\ll! \newdimen \uu! \newdimen\dd! \newdimen \rr! \newdimen
\temp!

\def\sq!#1#2#3#4#5{%
\ll!=#1 \uu!=#2 \dd!=#3 \rr!=#4
\setbox0=\hbox{%
 \temp!=\squaresize\advance\temp! by .5\uu!
 \rlap{\kern -.5\ll!
 \vbox{\hrule height \temp! width#1 depth .5\dd!}}%
 \temp!=\squaresize\advance\temp! by -.5\uu!
 \rlap{\raise\temp!
 \vbox{\hrule height #2 width \squaresize}}%
 \rlap{\raise -.5\dd!
 \vbox{\hrule height #3 width \squaresize}}%
 \temp!=\squaresize\advance\temp! by .5\uu!
 \rlap{\kern \squaresize \kern-.5\rr!
 \vbox{\hrule height \temp! width#4 depth .5\dd!}}%
 \rlap{\kern .5\squaresize\raise .5\squaresize
 \vbox to 0pt{\vss\hbox to 0pt{\hss $#5$\hss}\vss}}%
}
 \ht0=0pt \dp0=0pt \box0
}

\def\vsq!#1#2#3#4#5\endvsq!{\vbox to \squaresize{\hrule
width\squaresize height 0pt%
\vss\sq!{#1}{#2}{#3}{#4}{#5}}}

\newdimen \LL! \newdimen \UU! \newdimen \DD! \newdimen \RR!

\def\vvsq!{\futurelet\next\vvvsq!}
\def\vvvsq!{\relax
  \ifx     \next l\LL!=\Thickness \let\continue=\skipnexttoken!
  \else\ifx\next u\UU!=\Thickness \let\continue=\skipnexttoken!
  \else\ifx\next d\DD!=\Thickness \let\continue=\skipnexttoken!
  \else\ifx\next r\RR!=\Thickness \let\continue=\skipnexttoken!
  \else\def\continue{\vsq!\LL!\UU!\DD!\RR!}%
  \fi\fi\fi\fi
  \continue}

\def\skipnexttoken!#1{\vvsq!}

\def\place#1#2#3{\vbox to 0pt{\vss
\rlap{\kern#1\squaresize
  \raise#2\squaresize\hbox{$#3$}}
\vss}}

\catcode `!=12

\begin{document}
\maketitle


\begin{abstract}
In this paper, we answer a question posed by Kurt Johansson,
 to find a PDE for the joint distribution of the Airy Process.
 The latter is a continuous stationary process, describing the motion of the outermost particle of
 the
 Dyson Brownian motion, when the number of particles get large,
 with space and time appropriately rescaled.  The question
 reduces to an asymptotic analysis on the equation
 governing the joint probability of the eigenvalues of coupled Gaussian
 Hermitian matrices.

 The differential equations lead to the asymptotic
 behavior of the joint distribution and the correlation for the
  Airy process at
 different times $t_1$ and $t_2$, when
 $t_2-t_1 \rightarrow \iy$.

\end{abstract}

\tableofcontents

\section{Main result}

The Dyson Brownian
 motion (see \cite{Dyson})
 $$
 \Bigl(\lb_1(t),\ldots, \lb_n(t)\Bigr)\in \BR^n,$$
  with transition density $p(t,\mu,\lb)$
 satisfies the diffusion equation
\bean
  \frac{\pl p }{\pl t}
  &=&
 \frac{1}{2}
   \sum_1^n
  \frac{\pl}{\pl \lb_i} \Phi(\lb)
   \frac{\pl}{\pl \lb_i} \frac{1}{\Phi(\lb)}p,
  \\
  &=&
   \sum_1^n\left(\frac{1}{2}\frac{\pl^2}{\pl
\lb_i^2}-   \frac{\pl}{\pl \lb_i}\frac{\pl \log
\sqrt{\Phi(\lb)}}{\pl \lb_i}\right)p
  \eean
   with
 $$\Phi(\lb)= \Dt^2(\lb) \prod_1^n e^{-\lb_i^2/a}
 .$$
%
%

 It corresponds to the
 motion of the eigenvalues $(\lb_1(t),\ldots, \lb_n(t))$
  of an Hermitian matrix $B$,
  evolving according to the Ornstein-Uhlenbeck process
\be
\frac{\pl P}{\pl t}=
 \sum_{i,j=1}^{n^2} \Bigl( \frac{1}{4}(1+\delta_{ij})
 \frac{\pl^2 }{\pl B_{ij}^2}
  + \frac{1}{a^2}\frac{\pl }{\pl B_{ij}} (B_{ij}P)
  \Bigr)
  ~,\ee
  with transition density ($c:= e^{-t/a^2}$)
$$\hspace{-.5cm}
P(t,\bar B,B)=Z^{-1}\frac{1}{(1-c^2)^{n^2/2}}
 e^{ -  \frac{1}{a^2(1-c^2)}\Tr (B- c  \bar B)^2
    },
 $$
 The $B_{ij}$'s in (1.1) denote the $n^2$ free parameters in the
 Hermitian matrix $B$.
In the limit $t\rightarrow\iy$, this distribution tends
to the stationary distribution
%
  $$
  Z^{-1} e^{-\frac{1}{a^2}\Tr B^2} dB=
  Z^{-1}  \Dt^2(\lb) \prod_1^n e^{-\frac{\lb_i^2}{a^2}}
      ~ d\lb_i .
   $$
With
 ~this invariant measure as initial condition,
one
  finds for the joint distribution: 
 \be
  P(B(0)\in dB_1, B(t) \in dB_2)
  =
   Z^{-1}\frac{dB_1dB_2}{(1-c^2)^{n^2/2}}
 e^{ -  \frac{1}{a^2(1-c^2)} \Tr (
   B_1^2-
    2cB_1B_2 +
    B_2^2  )
    }\label{transition-probability}
 .\ee
%

  Setting $a=1$,
the {\bf Airy process} is defined by an appropriate
rescaling
 of the largest eigenvalue $\lb_n$ in the Dyson diffusion,
 \be
 A(t) =\lim_{n\rightarrow \iy}
  \sqrt{2} n^{1/6}
   \left( \lb_n(n^{-1/3}t)-\sqrt{2n} \right)
   \label{AiryP},\ee
   in the sense of convergence of distributions for a finite number of
   $t$'s. This process was introduced by
   Pr\"ahofer and Spohn \cite{Spohn} in the context of
   poly-nuclear growth models and
    further investigated by Johansson \cite{Johansson}. Pr\"ahofer and Spohn showed
 the Airy process
is a stationary process with continuous sample paths; thus the
probability
 $P(A(t)\leq u)$ is independent of $t$, and is given by the Tracy-Widom
 distribution
 \cite{TW},
  \be
 P(A(t)\leq u)=F_2(u):=
  \exp\left(-\int^{\iy}_{u}(\al-u)q^2(\al)
d\al \right),
 \ee
  with $q(\alpha)$
a solution of the {\bf Painlev\'e II} equation,
 \be
   q''=\al q+2q^3 ~~\mbox{with}~~q(\al)\cong
 \left\{\begin{array}{l}
 -\frac{
  e^{-\frac{2}{3}  {\al}^{\frac{3}{2}}}}{2\sqrt \pi \al^{1/4}}
\mbox{\,\,for\,\,}\al\nearrow \infty  \\  \\
  \sqrt{-\al/2}\mbox{\,\,for\,\,}\al\searrow -\infty.
  \end{array}\right.
  \ee

 At MSRI
(sept 02), Kurt Johansson posed the question, whether a
PDE can be found for the joint probability of the Airy
process; see \cite{Johansson}. The present paper answers
this question, which enables us to derive the
asymptotics of the large time correlation of the Airy
Process. We thank Kurt Johansson for introducing us to
this process.

\begin{theorem}
 Given $t_1<t_2$, the joint probability for the Airy process

$$
 h(t_2-t_1;~\frac{y+x}{2},\frac{y-x}{2}):=\log P\left(A(t_1) <\frac{y+x}{2}, A(t_2)<
  \frac{y-x}{2}\right),
 $$
%
 satisfies a non-linear PDE\footnote{in terms of the Wronskian
  $\{f(y),g(y)\}_y:=f'(y)g(y)-f(y)g'(y)$.} in $x,y$ and $t=t_2-t_1$,
\be
  2t   \frac{\pl^3 h}{\pl t \pl x \pl y}
 =
\left(t^2 \frac{\pl}{\pl x}- x  \frac{\pl}{\pl y}\right)
 \left(\frac{\pl^2 h}{\pl x^2}-
  \frac{\pl^2 h}{\pl y^2}\right)
+8 \left\{  \frac{\pl^2 h}{\pl x \pl y} ,
 \frac{\pl^2 h}{\pl y^2} \right\}_y,
\ee
 with initial condition
  $$
  \lim_{t\searrow 0}h\left(t;~\frac{y+x}{2},\frac{y-x}{2}
   \right)=\log  F_2
  \left(\min(\frac{y+x}{2},
  \frac{y-x}{2})\right)
  .$$

  \end{theorem}

  \vspace{.7cm}

  \noindent {\bf Conjecture}
   {\it For any fixed $t>0,~x \in \BR$, the conditional probability
   satisfies:
}%
\be \lim_{z\rightarrow \iy}P(A(t)\geq x+z \mid A(0)\leq
-z)
  =0.
\ee

\vspace{.7cm}

\noindent Accepting this conjecture, we prove:

\begin{theorem}
 For large $t$, the joint probability admits the
 asymptotic series
\be
 P(A(0) \leq u, A(t)\leq v)=
  F_2(u)F_2(v)
   +  \frac{F'_2(u)F'_2(v)}{t^2}
   + \frac{ \Phi(u,v)\!+\!\Phi(v,u)}{t^4}+
    O\left(\frac{1}{t^6}\right),
\label{jointprobability}\ee
with ($q(u)$ is the function (1.5))
\bean \Phi(u,v)&:=&
 F_2(u)F_2(v)
\left(\begin{array}{l}
  \frac{1}{4} \left({\displaystyle\int}_{\!\!\!u}^{\iy}  q^2d\al\right)^2
              \left({\displaystyle\int}_{\!\!\!v}^{\iy}  q^2 d\al\right)^2
   \\ \\+~q^2(u)  \left( \frac{1}{4}q^2(v) -\frac{1}{2}
   \bigl({\displaystyle\int}_{\!\!\!v}^{\iy}
          q^2 d\al \bigr)^2 \right) \\ \\
    +{\displaystyle\int}_{\!\!\!v}^{\iy}d\al\bigl( 2(v-\al)q^2
                          + q'^2
                          -q^4
                          \bigr)
                          {\displaystyle\int}_{\!\!\!u}^{\iy}
                          q^2d\al
         \end{array}\right).
    \eean
    Moreover, the covariance for large $t$ behaves as
    \be
    E(A(t)A(0))- E(A(t))E(A(0))=\frac{1}{t^2}+\frac{c}{t^4}+...
    ~~,\label{covariance}\ee
     where
$$
c:=2\int\!\!\!
 \int_{\BR^2}\Phi(u,v)du~dv.
$$

\end{theorem}

\remark The equation (1.6) for the probability
 $$
 h(t_2-t_1;u,v):=\log P(A(t_1) \le u, A(t_2)\le v), ~~~t=t_2-t_1,
 $$
 takes on the alternative form in the variables $u$ and
 $v$,
\bea
 t\frac{\pl}{\pl t}
  \Bigl( \frac{\pl^2}{\pl u^2}-
         \frac{\pl^2}{\pl v^2}  \Bigr)h
  &=&
  \frac{\pl^3h}{\pl u^2\pl v}
    \Bigl(  2 \frac{\pl^2h}{\pl v^2}
        +\frac{\pl^2h}{\pl u\pl v}-
         \frac{\pl^2h}{\pl u^2}+u-v-t^2 \Bigr)
        \no \\
   &&  -
     \frac{\pl^3 h}{\pl v^2\pl u}
    \Bigl(  2 \frac{\pl^2h}{\pl u^2}
        +\frac{\pl^2h}{\pl u\pl v}-
         \frac{\pl^2h}{\pl v^2}-u+v-t^2 \Bigr)
     \no \\
    && +
      \Bigl(\frac{\pl^3h}{\pl u^3}
             \frac{\pl}{\pl v}
             -\frac{\pl^3h}{\pl v^3}
             \frac{\pl}{\pl u}\Bigr)
             \Bigl(\frac{\pl}{\pl u} +\frac{\pl}{\pl v}
             \Bigr)h~,
  \label{PDE}\eea
 with initial condition
  $$
  \lim_{t\searrow 0}h(t;u,v)=\log  F_2(\min(u,v))
  .$$
This equation enjoys an obvious $u\leftrightarrow v$
 duality.


The proof of this theorem is based on a PDE, which was
obtained
 in {\cite{AvM1} for
the spectrum of coupled random matrices.

In \cite{Spohn}, Spohn and Pr\"ahofer pose the question
about the asymptotics of the covariant functions of
$A(t)$ and $A(0)$ for large $t$. Moreover, in a very
recent paper, Tracy and Widom \cite{TW2} expressed the
joint distribution, for several times $t_1,\ldots,t_m$,
as the exponential of a certain integral; its integrand
involves traces of matrices, which satisfy a coupled
system of non-linear ODE's. The quantities involved are
entirely different and their methods are
functional-theoretical; it remains unclear what the
connection is between the two results.

\section{The spectrum of coupled random matrices}

 Consider a product
ensemble $(M_1,M_2)\in {\cal H}_n^2 :={\cal
H}_n\times{\cal H}_n$ of $n\times n$ Hermitian matrices,
equipped with a Gaussian probability measure,
\begin{equation}
c_ndM_1dM_2\,e^{-\frac{1}{2}{\rm Tr}(M^2_1+M_2^2-2cM_1M_2)},
\label{probability}\end{equation} where $dM_1dM_2$ is Haar measure
on the product ${\cal H}_n^2$, with each $dM_i$,
\begin{equation}
dM_1=\Dt_n^2(x)\prod^n_1dx_idU_1\mbox{\,\,and\,\,}
dM_2=\Dt_n^2(y)\prod^n_1dy_idU_2 \label{2}
\end{equation}
decomposed into radial and angular parts.
 In \cite{AvM1}, we define differential operators ${\cal
A}_k$, $ {\cal B}_k$ of ``weight" $k$, which form a closed Lie
algebra, in terms of the coupling constant $c$, appearing in
(\ref{probability}), and the boundary of the set
\begin{equation}
 E=E_1\times E_2:=\cup^r_{i=1}[a_{2i-1},a_{2i}]\times
\cup^s_{i=1}[b_{2i-1},b_{2i}]\subset \BR^2.
 \end{equation}
Here we only need the first few ones:
 \begin{equation}
\begin{tabular}{ll}
${\cal
A}_1=\displaystyle{\frac{1}{c^2-1}\left(\sum^r_1\frac{\pl}{\pl
a_j}+c\sum^s_1\frac{\pl}{\pl b_j}\right)}$&${\cal
B}_1=\displaystyle{\frac{1}{1-c^2}\left(c\sum^r_1\frac{\pl}{\pl
a_j}+\sum^s_1\frac{\pl}{\pl b_j}\right)}$\\ ${\cal
A}_2=\displaystyle{\sum^r_{j=1}a_j^{}\frac{\pl}{\pl
a_j}-c\frac{\pl}{\pl c}}$&${\cal B}_2=\displaystyle{\sum^s_{j=1}
b_j^{}\frac{\pl}{\pl b_j}-c\frac{\pl}{\pl c}}.$
\end{tabular}
\end{equation}
 In \cite{AvM1}, we prove the following
 theorem:

\vspace{0.1cm}

\begin{theorem} 
  Given the joint distribution
   \begin{equation}
P_n(E):=P(\mbox{all}(M_1\mbox{-eigenvalues)}\in E_1,~
\mbox{all}(M_2\mbox{-eigenvalues})\in E_2),\label{statistics}
\end{equation}
 the function
  $F_n(c;a_1,\ldots,a_{2r},b_1,\ldots,b_{2s}): = \log P_n(E)$,
  satisfies the
    non-linear third-order partial
  differential equation\footnote{in terms of the Wronskian
  $\{f,g\}_X=Xf.g-f.Xg$, with regard to a first order
  differential operator $X$.}:
\be
 \left\{{\cal B}_2   {\cal A}_1 F_n  ~,~
  {\cal B}_1 {\cal A}_1 F_n +\frac{nc}{c^2-1} \right\}
_{ {\cal A}_1} ~-~ \left\{{\cal A}_2   {\cal B}_1 F_n
 ~,~ {\cal
A}_1 {\cal B}_1 F_n+\frac{nc}{c^2-1}\right\} _{ {\cal B}_1 }=0
 .
 \label{abequations}\ee
\end{theorem}


\begin{corollary} For $E=E_1\times E_2:=
 (-\infty,a]\times (-\infty,b]$,
  it is convenient to use the new variables
 $x:=-a+cb,~y:= -ac+b$. In these variables, the
  equation (\ref{abequations}) for
  \be f_n(c;x,y):=  \log P_n(E)=F_n\left(c; \frac{x-cy}{c^2-1},
  \frac{cx-y}{c^2-1}\right)
  \ee
  takes on the following form:
   \bea \lefteqn{
\frac{\pl}{\pl x}\left( \frac{(c^2-1)^2 \frac{\pl^2 f_n}{\pl x \pl
c}+2ncx-n(1+c^2)y}{(c^2-1) \frac{\pl^2 f_n}{\pl x
\pl y}+nc} \right) }\no\\
 &=&
 \frac{\pl}{\pl y}\left( \frac{(c^2-1)^2
\frac{\pl^2 f_n}{\pl y \pl c}+2ncy-n(1+c^2)x}{(c^2-1) \frac{\pl^2
 f_n}{\pl y \pl x}+nc} \right).
 \label{xyequations}\eea

\end{corollary}

\proof
 It is an immediate consequence of Theorem 2.1, upon observing
 the simple form of the differential operators
 ${\cal A}_1=\pl / \pl x$ and ${\cal B}_1=\pl / \pl y$,
 when expressed in terms of $x$ and $y$.\qed

\section{Proof of Theorem 1.1}

Taking into account the scaling in the
  definition of the Airy process (\ref{AiryP}),
  and using the Ornstein-Uhlenbeck transition probability
   (\ref{transition-probability}), we compute the probability
  (setting $c=e^{-n^{-1/3}t}$, with $t=t_2-t_1$)
\bea \lefteqn{P\left(\sqrt{2} n^{1/6}
   \bigl( \lb_n(n^{-1/3}t_1)-\sqrt{2n} \bigr)\leq u,
   \sqrt{2} n^{1/6}
   \bigl( \lb_n(n^{-1/3}t_2)-\sqrt{2n} \bigr)\leq v
   \right)
   }\no\\&&\no\\
&=& 
 \int\!\!\!\int_
  {\begin{array}{l}
   \mbox{all $B_1$-eigenvalues}
    \leq  \frac{1}{\sqrt{2}}
   (2 {n}^{1/2}+ {n^{-1/6}}  u) \\
  \mbox{all $B_2$-eigenvalues}
   \leq   \frac{1}{\sqrt{2}}
   (2 {n}^{1/2}+ {n^{-1/6}}   v)
    \end{array}  }
     \no\\
&&\hspace{5cm}
  Z^{-1}\frac{dB_1dB_2}{(1-c^2)^{n^2/2}}
 e^{ - \frac{1}{1-c^2}  \Tr (
   B_1^2  +  B_2^2-
    {2c} B_1B_2 )
    }
  \no\\
&=&
 \int\!\!\!\int_
  {\begin{array}{l}
   \mbox{all $M_1$-eigenvalues}
    \leq  \frac{2
 {n}^{1/2}+ {n^{-1/6}u}}{\sqrt{1-c^2}} \\
  \mbox{all $M_2$-eigenvalues}
   \leq   \frac{2
 {n}^{1/2}+ {n^{-1/6}v}}{\sqrt{1-c^2}}
    \end{array}  }
\no\\ &&\hspace{5cm}
Z'^{-1}  dM_1dM_2\,e^{-\frac{1}{2}{\rm Tr}(M^2_1+M_2^2-2cM_1M_2)},
 \no\\&&\label{2-matrixIntegral}\eea
 using the change of variables
 $$
 M_i=\frac{  B_i}{\sqrt{(1-c^2)/2}}.
 $$
The integral (\ref{2-matrixIntegral}) turns out to coincide with
the
 statistics (\ref{statistics}) of the coupled matrix model with the
  Gaussian distribution (\ref{probability}).
 Therefore, we set \be
 a=\frac{2
 {n}^{1/2}+ {n^{-1/6}u}}{\sqrt{1-c^2}}~~\mbox{and}~~
 b=\frac{2
 {n}^{1/2}+ {n^{-1/6}v}}{\sqrt{1-c^2}};
 \ee
 in formula (\ref{xyequations}), via $x$ and $y$. Then, setting $k=n^{1/6}$, we now express $x$ and $y$ in terms of
$u$, $v$ and $c=e^{-n^{-1/3}t}=e^{-t/k^2}$, using the change of
variables in Corollary 2.2,
\bea
x(c(t),u,v)&=&-a+bc=-\frac{1}{\sqrt{1-c^2}}\left(\left(2k^3+\frac{u}{k}\right)
-\left(2k^3+\frac{v}{k}\right)c\right)\no\\
\no\\
y(c(t),u,v)&=&-ac+b=-\frac{1}{\sqrt{1-c^2}}\left(\left(2k^3+\frac{u}{k}\right)c
-\left(2k^3+\frac{v}{k}\right)\right)
 \no\\ \label{directmap}\eea
 with inverse given by
\bea
u&=&\frac{1}{\sqrt{1-c^2}}(k(cy-x)-2k^4\sqrt{1-c^2})\no\\
 \no\\
v&=&\frac{1}{\sqrt{1-c^2}}(k(y-cx)-2k^4\sqrt{1-c^2}),
\quad\mbox{with}~c=e^{-t/k^2}.\no\\
 \label{inversemap} \eea
 So, $f_n(c;x,y)$ satisfies equations (\ref{xyequations}), which
 written out, has the form
\bea \lefteqn{\frac{\pl^3f}{\pl x^2\pl
y}\left((c^2-1)^2\frac{\pl^2 f}{\pl x\pl c}
+2cnx-(c^2+1)ny\right)}\no\\
 \no\\
& &-\frac{\pl^3  f}{\pl x\pl y^2}\left((c^2-1)^2\frac{\pl^2 f}{\pl
y\pl c}+2cny-
(c^2+1)nx\right)\no\\
\no\\
& &+(c^2-1)\left((c^2-1)\frac{\pl^2 f}{\pl x\pl
y}+nc\right)\left(\frac{\pl^3  f}{\pl y^2\pl c}-\frac{\pl^3 f}{\pl
x^2\pl c}\right)=0. \label{f-equation}
 \eea
 Using (\ref{directmap}) and definition (2.7), one defines, at a first stage, a new
function \bean
  g_n(c(t);u
      ,v
      )&:=&f_n(c(t);x(c(t);u,v),y(c(t);u,v))  \\
 &=&  F_n\left(c(t);\frac{x(c;u,v)-cy(c;u,v)}{c^2-1},
 \frac{cx(c;u,v)-y(c;u,v)}{c^2-1}\right),
 \eean
 one expresses the $(x,y,c)$-partials of $f$ in terms
 of $(u,v,c)$-partials of $g$, using the inverse map
 (\ref{inversemap}), e.g.,
 $$
 \frac{\pl^3f}{\pl x^2\pl  y}=
  \left(\frac{k}{\sqrt{1-c^2}}\right)^3
  \left(
   c( \frac{\pl^3g}{\pl u^3}\! +\! c\frac{\pl^3g}{\pl v^3})
    +(2c^2\!+\!1)\frac{\pl^3g}{\pl u^2\pl v}
    +c(c^2\!+\!2)\frac{\pl^3g}{\pl u\pl v^2}
    \right),
 $$
 etc... .One substitutes these expressions in
 the differential equation (\ref{f-equation}) for $f(c;x,y)$,
 yielding a differential equation in $g(c;u,v)$, with coefficients
 depending on $u,v,c$ and $k$. After division by $k^4$,
 this differential equation
  is a quadratic polynomial in $k^4$,

\bea
 \lefteqn{k^{8}(c-1)^2~2 \left(\frac{\pl^3 g}{\pl
   u^2\pl v}-
\frac{\pl^3 g}{\pl u\pl v^2}\right)}\no\\
&&
 \no\\
&&
 +k^4\left\{\begin{array}{l}
%
%
%
\displaystyle{+\frac{\pl^3 g}{\pl u^2\pl
v}\left(2(c^2-c+1)\frac{\pl^2  g}{\pl
u^2}-4c\frac{\pl^2  g}{\pl v^2}+v(c^2+1)-2cu\right)}\\
\\
-\displaystyle{\frac{\pl^3  g}{\pl u\pl
v^2}\left(2(c^2-c+1)\frac{\pl^2  g}{\pl
v^2}-4c\frac{\pl^2  g}{\pl u^2}+u(c^2+1)-2cv\right)}\\
\\
  -2c\displaystyle{\left(\frac{\pl^3 g}{\pl u^3}\frac{\pl^2  g}{\pl
v^2}-\frac{\pl^3  g}{\pl v^3}\frac{\pl^2  g}{\pl
u^2}+\Bigl(\frac{\pl^3  g}{\pl u^2\pl v}-\frac{\pl^3 g}{\pl
u\pl v^2}\Bigr)\frac{\pl^2 g}{\pl u\pl v}\right)}\\
\\
- 2(c^2-c+1)\displaystyle{\frac{\pl^2 g}{\pl u\pl
v}\left(\frac{\pl^3 g}{\pl
u^3}-\frac{\pl^3 g}{\pl v^3}\right)}
 \\
\\
 +c(c^2-1)\displaystyle{\left(\frac{\pl^3 g}{\pl u^2\pl
c}-\frac{\pl^3 g}{\pl v^2\pl c}\right)}
\end{array}
\right\} +\ldots.\no\\&&
 \label{k-expansion}
  \eea
 Then taking into account the fact that
\be c(t)=e^{-t/k^2}=1-\frac{t}{k^2}+{
O}\left(\frac{1}{k^4}\right),\quad\frac{\pl}{\pl
c}=-\frac{k^2}{c} \frac{\pl}{\pl t},
 \label{c-expansion}
 \ee
 the leading term has order $k^4$ for large $k$
 (noting that $k^4(c-1)^2=t^2+O(1/k^2)$),
 At a second stage, defining
$$
h(t;u,v):=\lim_{k\rg\iy}g(e^{-t/k^2};u,v),
$$
and using the expansion (\ref{c-expansion}) for $c(t)$
and
 the partial ${\pl}/{\pl c}=-({k^2}/{c})~
  {\pl}/{\pl t}$,
 the leading
term in the expression (\ref{k-expansion}) has order
$k^4$; no contribution comes from the $k^0$-term. This
leading term (multiplied with $-1/2$) must therefore
vanish, leading first to equation (\ref{PDE}) and then
setting $x=u-v,~y=u+v$ in that equation, to equation
(1.6) given in Theorem 1.\qed
\newpage

\section{Proof of Theorem 1.2}


Before giving the proof of Theorem 1.2, we remind the
reader of the conjecture stated in section 1: for any
fixed $t>0,~x \in \BR$, the conditional probability
   satisfies:
\be \lim_{z\rightarrow \iy}P(A(t)\geq x+z \mid A(0)\leq
-z)
  =0.
\ee
 Under this assumption, we prove the following:

\begin{lemma}
 Considering the series for the probability, for large $t$,
   \be P(A(0)\leq u,A(t)\leq
v)=F_2(u)F_2(v)\left(1+\sum_{i\geq
1}\frac{f_i(u,v)}{t^i}\right),
 \ee
  the coefficients
$f_i(u,v)$ have the property
 \be
\lim_{u\rg\iy}f_i(u,v)=
 \lim_{u\rg\iy}f_i(v,u)=0,~\mbox{for fixed $v\in \BR$}
\ee and \be \lim_{z\rg\iy}f_i(-z,z+x)=0,~~
 ~~\mbox{for fixed $x\in \BR$}.
  \ee

\end{lemma}

\proof First observe by (5.1) in section 5 that the Airy
kernel becomes diagonal when $t\rightarrow \iy$. Then
the Airy process decouples at $\iy$, and, using the
stationarity, one is lead to
$$
\lim_{t\rightarrow \iy} P(A(0)\leq u,A(t)\leq
v)=P(A(0)\leq u)P(A(0)\leq v)=F_2(u)F_2(v).
$$
 The next terms follow from the PDE (1.6),
  although it is more convenient here to use
   the form (1.10) of the equation.
%
%
Considering the following conditional probability,
 \bean
\lefteqn{P(A(t)\leq v\mid A(0)\leq u)}\\
\\
&=&\frac{P(A(0)\leq u,A(t)\leq v)}{P(A(0)\leq u)}\\
\\
&=&F_2(v)\left(1+\sum_{i\geq
1}\frac{f_i(u,v)}{t^i}\right),
 \eean
 and setting
  $$
 v=z+x,\qquad u=-z ,
  $$
%
%
 we have for all $t$, since $\lim_{z\rg \iy}
 F_2(z+x)=1$, and by (4.1) that
 \bean
1&=&\lim_{z\rg\iy}P(A(t)\leq z+x\mid A(0)\leq -z)\\
\\
&=& 1+\sum_{i\geq
1}\frac{\displaystyle{\lim_{z\rg\iy}}f_i(-z,z+x)}{t^i}
,\eean
implying that
 $$
\lim_{z\rg\iy}f_i(-z,z+x)=0,\qquad\mbox{ for all } i\geq
1.
 $$

 Similarly, letting $v\rg\iy$, we have \bean
1=\lim_{v\rg\iy}P(A(t)\leq v\mid A(0)\leq
u)&=&\lim_{v\rg\iy}\left[F_2(v)\left(1+\sum_{i\geq
1}\frac{f_i(u,v)}{t^i}\right)\right]\\
\\
&=&1+\sum_{i\geq
1}\frac{\displaystyle{\lim_{v\rg\iy}}f_i(u,v)}{t^i}.
\eean Hence
  \be
   \lim_{v\rg\iy}f_i(u,v)=0
   \ee
    and, considering the same argument for the conditional probability
    $P(A(0)\leq u\mid A(t)\leq v)$,
$$
\lim_{u\rg\iy}f_i(u,v)=0,
$$
ending the proof of Lemma 4.1.\qed

\vspace{1cm}

\noindent{\it Proof of Theorem 1.2:\/}
Putting the log of the expansion (4.2)
  \bea
h(t;u,v)&=&\log P(A(0)\leq u,A(t)< v)\nonumber\\
&=&\log F_2(u)+\log F_2(v)+\sum_{i\geq
1}\frac{h_i(u,v)}{t^i}  \no\\
 &=&\log F_2(u)+\log F_2(v)+\frac{f_1(u,v)}{t}+
  \frac{f_2(u,v)-f_1^2(u,v)/2}{t^2}+\ldots
,\no\\
  \eea
in the equation (\ref{PDE}), leads to:

\medbreak

 \noindent{\bf (i)} a leading term of order
$t$, given by
  \be  \LR
 h_1=0,
 \ee
 where
 \be
 \LR:=
   \left(\frac{\pl}{\pl u}-\frac{\pl}{\pl
v}\right)\frac{\pl^2}{\pl u\pl v}.
 \ee
  The most general solution to (4.7) is given by
$$
h_1(u,v)=r_1(u)+r_3(v)+r_2(u+v),
$$
with arbitrary functions $r_1,r_2,r_3$.
Hence,
$$
P(A(0)\leq u,A(t)\leq v)
=F_2(u)F_2(v)\left(1+\frac{h_1(u,v)}{t}+...\right)
$$
with $h_1(u,v)=f_1(u,v)$ as in (4.2). Applying Lemma
4.1,
$$
r_1(u)+r_3(\iy)+r_2(\iy)=0,\quad\mbox{for all~}u\in\BR,
$$
implying
$$
r_1(u)=\mbox{constant}=r_1(\iy),
$$
and similarly
$$
r_3(u)=\mbox{constant}=r_3(\iy).
$$

Therefore, without loss of generality, we may absorb the
constants $r_1(\iy)$ and $r_3(\iy)$ in the definition of
$r_2(u+v)$. Hence, from (4.6),
$$
f_1(u,v)=h_1(u,v)=r_2(u+v)
$$
using (4.5),
$$
0=\lim_{z\rg\iy}f_1(-z,z+x)=r_2(x)
$$
implying that the $h_1(u,v)$-term in the series (4.6)
vanishes.

\medbreak

\noindent{\bf (ii)} One computes that the term
$h_2(u,v)$ in the expansion (4.6) of $h(t;u,v)$
satisfies
 \be
\LR
  h_2=\frac{\pl^3g}{\pl
u^3}\frac{\pl^2g}{\pl v^2}- \frac{\pl^3g}{\pl
v^3}\frac{\pl^2g}{\pl u^2}, \mbox{~with~}g(u):=\log
F_2(u).
 \ee
This is the term of order $t^0$, by putting the series
(4.6) in the equation (\ref{PDE}). The most general
solution
to (4.9) is
$$
h_2(u,v)=g'(u)g'(v)+r_1(u)+r_3(v)+r_2(u+v).
$$
Then \bean
P(A(0)\leq u,A(t)\leq v)&=&e^{h(t,u,v)}\\
\\
&=&F_2(u)F_2(v)e^{\displaystyle{\sum_{i\geq 2}\frac{h_i(u,v)}{t^i}}}\\
\\
&=&F_2(u)F_2(v)\left(1+\frac{h_2(u,v)}{t^2}+...\right).
\eean
 In view of the explicit formula for
  the distribution $F_2$
 and the behavior (1.5) of $q(\al)$ for
$\al\nearrow\iy$, we have that
  \bean
\lim_{u\rg\iy}g'(u)&=&\lim_{u\rg\iy}
 (\log F_2(u))^{\prime}\\
&=&\lim_{u\rg\iy}\int_u^{\iy}q^2(\al)d\al =0. \eean
Hence
$$
0=\lim_{u\rg\iy}f_2(u,v)=\lim_{u\rg\iy}h_2(u,v)
 =r_1(\iy)+r_3(v)+r_2(\iy),
$$
showing $r_1$ and similarly $r_3$ are constants.
Therefore, by absorbing $r_1(\iy)$ and $r_3(\iy)$ into
$r_2(u+v)$, we have
$$
f_2(u,v)=h_2(u,v)=g'(u)g'(v)+r_2(u+v).
$$
Again, by the behavior of $q(x)$ at $+\iy$ and $-\iy$,
$$
g'(-z)g'(z+x)=\int^{\iy}_{-z}q^2(\al)d\al\int_{z+x}^{\iy}q^2(\al)d\al\leq
cz^{3/2}e^{-2z/3}.
$$
Hence
$$
0=\lim_{z\rg\iy}f_2(-z,z+x)=r_2(x)
$$
and so
$$
f_2(u,v)=h_2(u,v)=g'(u)g'(v),
$$
yielding the $1/t^2$ term in the series (4.6).

\medbreak

\noindent{\bf (iii)} Next, setting
 \bea
h(t;u,v)&=&\log P(A(0)\leq u,A(t)\leq v)\nonumber\\
\nonumber\\
&=&g(u)+g(v)+\frac{g'(u)g'(v)}{t^2}+
\frac{h_3(u,v)}{t^3}+...
 \eea
  in the equation (\ref{PDE}), we find for the $t^{-1}$ term:
 $$
\LR
 h_3=0.
 $$
As in (4.7), its most general solution is given by
$$
h_3(u,v)=r_1(u)+r_3(v)+r_2(u+v).
$$
By exponentiation of (4.6), we find
  \bean
   P(A(0)\leq u,A(t)\leq v)
    =F_2(u)F_2(v)\Biggl(1+\frac{g'(u)g'(v)}{t^2}\\
 \hspace*{3cm} +\frac{r_1(u)+r_3(v)+r_2(u+v)}{t^3}+...
 \Biggr).
 \eean
 The precise same arguments lead to
  $h_3(u,v)=0$.

\medbreak

 \noindent{\bf (iv)} So, at the next stage, we have
  \be
  h(t;u,v)=g(u)+g(v)+\frac{g'(u)g'(v)}{t^2}+
\frac{h_4(u,v)}{t^4}+\ldots
  \ee
  with
  \be
  f_4(u,v)=h_4(u,v)+\frac{1}{2} h^2_2(u,v)
  =h_4(u,v)+\frac{1}{2} g'(u)^2g'(v)^2.
  \ee
  Setting the series (4.11) in the equation (\ref{PDE}), we find
  for the $t^{-2}$ term:
   \bea
   \LR h_4
  &=&
  2\left(\frac{\pl^3g}{\pl
u^3}\left(\frac{\pl^2g}{\pl
v^2}\right)^2-\frac{\pl^3g}{\pl
v^3}\left(\frac{\pl^2g}{\pl
u^2}\right)^2\right)
 +\frac{\pl^3g}{\pl
u^3}\frac{\pl^3g}{\pl v^3}\left(\frac{\pl g}{\pl u}
 -\frac{\pl g}{\pl v}\right)\nonumber\\
& &+\frac{1}{2}\left(\frac{\pl^4g}{\pl
u^4}\frac{\pl}{\pl v}\left(\frac{\pl g}{\pl
v}\right)^2-\frac{\pl^4g}{\pl  v^4}\frac{\pl}{\pl
u}\left(\frac{\pl g}{\pl u}\right)^2\right)\nonumber\\
& &+\left(\frac{\pl^3g}{\pl
u^3}\frac{\pl^2g}{\pl v^2}+\frac{\pl^3g}{\pl  v^3}
\frac{\pl^2g}{\pl u^2}\right)(u-v)
 +2\left(\frac{\pl^3g}{\pl u^3}\frac{\pl g}{\pl v}
-\frac{\pl^3g}{\pl  v^3} \frac{\pl g}{\pl u}\right)
 \nonumber\\\no\\
&=&
 2\Big(2q(u)q'(u)(q(v)q'(v)+1)-q(u)q''(u)q^2(v)
   -(q'(u))^2q^2(v)\Big)\int_v^{\iy}q^2  \no \\
& &
 +~2q(u)\Big(q(u)q'(v)q''(v)+q'(u)q(v)q''(v)
  -2q(u)q^3(v)q'(v)\Big)\no\\
& &-\mbox{~same with~}u\leftrightarrow v.
 \eea
  This latter is an expression in $q(u)$, $q(v)$ and its
derivatives and in
$\displaystyle{\int_u^{\iy}}q^2(\al)d\al$ and
$\displaystyle{\int_v^{\iy}}q^2(\al)d\al$. It is
obtained by substituting in the previous expression
$$
g(u)=\int_u^{\iy}(u-\al)q^2(\al)d\al
$$
and the Painlev\'e II differential equation for $q(u)$,
$$
u~q(u)=q''(u)-2q(u)^3,
$$
in order to eliminate the explicit appearance of $u$ and
$v$.

\newpage

Now introducing\footnote{Note
$$
g'(u)=\int_u^{\iy}q^2(\al)d\al,~~~~g''(u)=-q^2(u),
$$
and
$$
g'_1(u)=\int_u^{\iy}q^{\prime
2}(\al)d\al,~~~~~g'_2(u)=\int_u^{\iy}q^4(\al)d\al.
$$
}
\bean
g(u)&=&\int_u^{\iy}(u-\al)q^2(\al)d\al\\
g_1(u)&=&\int_u^{\iy}(u-\al)q^{\prime 2}(\al)d\al\\
g_2(u)&=&\int_u^{\iy}(u-\al)q^4(\al)d\al, \eean the most
general solution to equation (4.13) is given, modulo the
null-space of $\LR$, by

\bea
  h_4(u,v)&=&{\frac{1}{2}\Big(g''(u)g'(v)^2+
   g''(v)g'(u)^2+g''(u)g''(v)\Big)}\no\\
& &+~g'(u)\Big(2g(v)+g'_1(v)-g'_2(v)\Big)  \no\\
 & &   +~g'(v)\Big(2g(u)+g'_1(u)-g'_2(u)\Big)\no\\
&=&q^2(u)\left(\frac{q^2(v)}{4}-\frac{1}{2}
 \left(\int_v^{\iy}q^2(\al)d\al\right)^2\right)\no\\
&
&+\int^{\iy}_uq^2(\al)d\al\int_v^{\iy}\Big(2(v-\al)q^2(\al)+q^{\prime
2}(\al)-q^4(\al)\Big)d\al  \no\\
& &  \no\\
& & +\mbox{~same with~}u\leftrightarrow v.
   \eea
 This form, together with (4.12), implies for the
 function $f_4(u,v)$:
\bean f_4(u,v)&=&h_4(u,v)+\frac{1}{2} g'(u)^2g'(v)^2
                       +r_1(u) +r_3(v)+r_2(u+v)  \\
&=&\sum_i a_i(u)b_i(v)+r_1(u)+r_3(v)+r_2(u+v).
  \eean
Using the asymptotics of $q(u)$, one finds

\bean
a_i(u),b_i(u)&\leq& c~e^{-u}\qquad u\rg\iy,\\
&\leq&c |u|^3\qquad u\rg -\iy, \eean and so, by the same
argument,
$$
r_1(u)=r_2(u)=r_3(u)=0.
$$
Therefore, we have
 $$ f_4(u,v)=h_4(u,v)+\frac{1}{2} g'(u)^2g'(v)^2.$$
with $h_4(u,v)$ as in (4.14), thus yielding the formula
(\ref{jointprobability}).


Finally, to prove formula (\ref{covariance}), we
compute, after integration by parts,
 \bean
E\Big(A(0)A(t)\Big)&=&\int\!\!\!\int_{\BR^2}uv
 \frac{\pl^2}{\pl u~\pl v}P(A(0)\leq u,A(t)\leq v)du~dv\\
&=&\int^{\iy}_{-\iy}uF'_2(u)du\int^{\iy}_{-\iy}vF'_2(v)dv\\
& &+~\frac{1}{t^2}
 \int^{\iy}_{-\iy}F'_2(u)du\int^{\iy}_{-\iy}F'_2(v)dv\\
& &+~\frac{1}{t^4}\int\!\!\!\int_{\BR^2}
 \Big(\Phi(u,v)+\Phi(v,u)\Big)du~dv\\
& &+~\mbox{ O}\left(\frac{1}{t^6}\right)\\
&=&\Bigl(E\Big(A(0)\Big)\Bigr)^2+\frac{1}{t^2}
 +\frac{c}{t^4}+\mbox{ O}\left(\frac{1}{t^6}\right),
  \eean
  where
$$
c:=\int\!\!\!\int_{\BR^2}\Big(\Phi(u,v)+\Phi(v,u)\Big)
  du~dv
=2\int\!\!\!\int_{\BR^2}  \Phi(u,v)dudv,$$ thus ending
the proof of Theorem 1.2. \qed


\section{The extended Airy kernel}
 The joint probabilities for the Airy process can also
 be expressed in terms of the Fredholm determinant of a matrix
 kernel, the so-called extended Airy kernel
  (considered in \cite{Johansson}, \cite{Forrester} and
 \cite{Spohn}), namely

  \bean
 \lefteqn{ P(A(t_1) <u_1,\ldots, A(t_m)<u_m)} \\
  &=&\left.\det\left( I-z (\hat K_{ij})_{1\leq i,j\leq m}
   \right)\right|_{z=1}
  \\
    &=& \left.1+ \sum_{k=1}^{\iy}\frac{(-z)^k}{k!}
  \sum_{1\leq i_1\leq \ldots \leq i_k\leq m}\int_{\BR^k}
   \det \Bigl(  \hat K_{i_ri_s}(x_r,x_s)\Bigr)_{1\leq r,s\leq k}
  dx_1\ldots dx_k \right|_{z=1}
  \eean
  with
  \be
   \hat K_{ij}(x,y):= \chi_{[u_i,\iy)}(x)
                K_{ij}(x,y)
                \chi_{[u_j,\iy)}(y)
    \ee
    and
    $$
     K_{ij}(x,y):=\left\{
       \begin{array}{l}
         \int_0^{\iy} e^{-z (t_i-t_j)} \mbox{Ai}(x+z)
                                        \mbox{Ai}(y+z)dz,
         ~~~~~~\mbox{if} ~t_i\geq t_j
           \\ \\
         - \int^0_{-\iy} e^{-z (t_i-t_j)} \mbox{Ai}(x+z)
                                        \mbox{Ai}(y+z)dz,
                                       ~~ \mbox{if} ~t_i < t_j
        ~,\end{array}\right.
 $$
Ai$(x)$ being the Airy function. So the Fredholm
determinant above is also a solution of the PDE of
Theorem 2.1 for $m=2$.


\section{Appendix: remarks about the conjecture}

Consider the Dyson Brownian motion $(\lb_1(t),\ldots,
\lb_n(t))$ and the corresponding Ornstein-Uhlenbeck
process on the matrix $B$. Then, using the change of
variables
 $$
 M_i=\frac{  B_i}{\sqrt{(1-c^2)/2}},
 $$ and further
 $M_2\mapsto M:=M_2-cM_1$  in the
$M_2$-integrals below and
 noting that $\max(\mbox{sp}~M_1) \leq -z$ and
       $\max(\mbox{sp}~M_2) \geq a$ imply
       $\max(\mbox{sp}~(M_2-cM_1) \geq a +cz $, we have
       for the conditional probability, the following
       inequality:
\bean
 \lefteqn{P\left( \lb_n(t)\geq a ~\big|~ \lb_n(0)\leq -z \right)
}\\ &&\\
 & =&
 \frac{
 {\displaystyle\int}_{\max(\mbox{\footnotesize sp}~ M_1)
                                             \leq -z}dM_1
    e^{-\frac{1}{2} (1-c^2)\Tr M_1^2}
{\displaystyle\int}_{\max(\mbox{\footnotesize sp}~M_2)
                                            \geq a}dM_2
     e^{-\frac{1}{2} \Tr (M_2-cM_1)^2}}
      {{\displaystyle\int}_{\max(\mbox{\footnotesize sp}~ M_1)\leq -z}dM_1
    e^{-\frac{1}{2} (1-c^2)\Tr M_1^2}
     {\displaystyle\int}_{M_2\in {\cal H}_n}dM_2
     e^{-\frac{1}{2} \Tr (M_2-cM_1)^2}}   \\ &&\\
  &\leq&
 \frac{{\displaystyle\int}_{\max(\mbox{\footnotesize sp}~ M_1)\leq -z}dM_1
    ~e^{-\frac{1}{2} (1-c^2)\Tr M_1^2}
     {\displaystyle\int}_{\max(\mbox{\footnotesize sp}~
                     M)\geq a+cz}dM
     ~e^{-\frac{1}{2} \Tr M^2}}
      {{\displaystyle\int}_{\max(\mbox{\footnotesize sp}~ M_1)\leq -z}dM_1
    ~e^{-\frac{1}{2} (1-c^2)\Tr M_1^2}
     {\displaystyle\int}_{M\in {\cal H}_n}dM
     ~e^{-\frac{1}{2} \Tr M^2}} \\
     &&
     \\
     &=&
     P(\lb_n(t)\geq a+cz),
 \eean
implying
 $$
 \lim_{z\rightarrow \iy}
 P\left( \lb_n(t)\geq a ~\big|~ \lb_n(0)\leq -z
 \right)=0
,
 $$
 and a fortiori,
  $$
 \lim_{z\rightarrow \iy}
 P\left( \lb_n(t)\geq x+z ~\big|~ \lb_n(0)\leq -z
 \right)=0
.
 $$
 It is unclear why this limit remains valid when $n\rightarrow \iy$,
 using the Airy scaling (1.3). But the extended
 Airy kernel (5.1) seems to indicate the conjecture is valid.


\newpage

\end{document}